\documentclass[12pt]{article}
\newcommand{\ol}{\setlength{\itemsep}{0pt.}\begin{enumerate}}
\newcommand{\eol}{\end{enumerate}\setlength{\itemsep}{-\parsep}}

\newtheorem{thm}{Theorem}
\newtheorem{lm}{Lemma}
\newtheorem{pr}{Proposition}

\newtheorem{re}{Remark}

\newtheorem{problem}{Problem}

\usepackage{amsmath, amsfonts}

\newcommand{\e}{\varepsilon}
\newcommand{\N}{\mathbb{N}}
\newcommand{\R}{\mathbb{R}}
\newcommand{\1}{ \mathbf{1} }

\def\calU{{\cal U}}

\def\pf{ \noindent {\bf Proof: \  }}

\newcommand{\qed}{\hfill\vrule height6pt
width6pt depth0pt}

\def\endpf{\qed \medskip}

\renewcommand{\qed}{\hfill\vrule height6pt  width6pt depth0pt}

\def \an{\{ a_n \}_{n=1}^{\infty}}
 \def \xn{\{ x_n \}_{n=1}^{\infty}}

\def \xnk{\{ x_{n(k)} \}_{k=1}^{\infty}}

\def \tn{\{ t_n \}_{n=1}^{\infty}}

\def \Txn{\{ Tx_n \}_{n=1}^{\infty}}

\def \wn{\{ w_n \}_{n=1}^{\infty}}

\def \vn{\{ v_n \}_{n=1}^{\infty}}

\def \yn{\{ y_n \}_{n=1}^{\infty}}

 \def \En{\{ E_n \}_{n=1}^{\infty}}

\def \Enxn{\{ \1_{E_n} x_n \}_{n=1}^{\infty}}

\def\ep{\varepsilon}

\def\ss#1{_{\lower3pt\hbox{$\scriptstyle #1$}}}

\title{{Multiplication operators on $L(L_p)$ and $\ell_p$-strictly singular operators}}
\author{William B. Johnson\thanks{Supported in part by NSF
 DMS-0503688 and U.S.-Israel Binational Science  Foundation.}\thanks{Corresponding author.}
\ \ Gideon Schechtman\thanks{Supported in part by Israel Science
Foundation and  U.S.-Israel Binational Science Foundation;
participant, NSF Workshop in  Analysis and Probability, Texas A\&M
University.} }
\begin{document}
\maketitle


\begin{abstract}
A classification of  weakly compact multiplication operators on
$L(L_p)$, $1<p<\infty$, is given.  This answers a question raised
by Saksman and Tylli in 1992.  The classification involves the
concept of $\ell_p$-strictly singular operators, and we also
investigate the structure of general $\ell_p$-strictly singular
operators on $L_p$. The main result is that if an operator $T$ on
$L_p$, $1<p<2$, is $\ell_p$-strictly singular and $T_{|X}$ is an
isomorphism for some subspace $X$ of $L_p$, then $X$ embeds into
$L_r$ for all $r<2$, but $X$ need not be isomorphic to a Hilbert
space.

It is also shown that if $T $ is convolution by a biased coin on
$L_p$ of the Cantor group, $1\le p <2$, and $T_{|X}$ is an
isomorphism for some reflexive subspace $X$ of $L_p$, then $X$ is
isomorphic to a Hilbert space.  The case $p=1$ answers a question
asked by Rosenthal in 1976.

\end{abstract}

\section{Introduction}

\ \ \ \ Given (always bounded, linear)  operators $A$, $B$ on a
Banach space $X$, define $L_A$, $R_B$ on $L(X)$ (the space of
bounded linear operators on $X$) by $L_AT=AT$, $R_BT=TB$.
Operators of the form $L_A R_B$ on $L(X)$ are called {\em
multiplication operators}.  The beginning point of this paper is a
problem raised in 1992 by E. Saksman and H.-O. Tylli \cite {st1}
(see also \cite[Problem 2.8]{st2}):

Characterize the multiplication operators on $L(L_p)$, $1<p\not= 2<\infty$, which are
weakly compact.

Here $L_p$ is $L_p(0,1)$ or, equivalently, $L_p(\mu)$ for any purely non-atomic separable
probability  $\mu$.

In Theorem \ref{multiplication}
 we answer the Saksman-Tylli question.  The characterization is
rather simple but gives rise to   questions about operators on
$L_p$, some of which were asked by Tylli. First we set some
terminology. Given an operator $T:X\to Y$ and a Banach space $Z$,
say that $T$ is $Z$-strictly singular provided there is no
subspace $Z_0$ of $X$  which is isomorphic to $Z$ for which
$T_{|{Z_0}}$ is an isomorphism.  An operator $S:Z\to W$ factors
through an operator $T:X\to Y$  provided there are operators
$A:Z\to X$ and $B:Y\to W$ so that $S=BTA$.  If $S$ factors through
the identity operator on $X$, we say that $S$ factors through $X$.

If $T$ is  an operator on $L_p $, $1<p<\infty$, then $T$ is
$\ell_p$-strictly singular (respectively, $\ell_2$-strictly
singular)  if and only if $I_{\ell_p}$ (respectively,
$I_{\ell_2}$) does not factor through $T$.  This is true because
every subspace of $L_p$ which is isomorphic to ${\ell_p}$
(respectively, ${\ell_2}$) has a subspace which is still
isomorphic to ${\ell_p}$ (respectively, ${\ell_2}$)  and is
complemented in $L_p$.   Actually, a stronger fact is true: if
$\xn$ is a sequence in $L_p$ which is equivalent to the unit
vector basis for either $\ell_p$ or $\ell_2$, then $\xn$ has a
subsequence which spans a complemented subspace of $L_p$. For
$p>2$, an even stronger theorem was proved by Kadec-Pe\l czy\'nski
  \cite{kp}.  When $1<p<2$ and $\xn$ is a sequence in
$L_p$ which is equivalent to the unit vector basis for $\ell_2$,
one takes $\yn$ in $L_{p'}$ which are uniformly bounded and
biorthogonal to $\xn$.  By passing to a subsequence    which is
weakly convergent and subtracting the limit from each $y_n$, one
may assume that $y_n \to 0$ weakly and hence, by the Kadec--Pe\l
czy\'nski dichotomy \cite{kp}, has a subsequence that is
equivalent to the unit vector basis of $\ell_2$ (since it is
clearly impossible that $\yn$ has a subsequence equivalent to the
unit vector basis of $\ell_{p'}$). This implies that that the
corresponding subsequence of $\xn$ spans a complemented subspace
of $L_p$. (Pe\l czy\'nski showed this argument, or something
similar, to one of the authors many years ago, and a closely related result was proved in \cite{pelros}.) Finally, when
$1<p<2$ and $\xn$ is a sequence in $L_p$ which is equivalent to
the unit vector basis for $\ell_p$, see the comments after the
statement of Lemma \ref{ellpsequence}.

Notice that the comments in the preceding paragraph yield that an operator on $L_p $,
$1<p<\infty$, is
$\ell_p$-strictly singular (respectively, $\ell_2$-strictly singular)  if and only if
$T^*$ is $\ell_{p'}$-strictly singular (respectively, $\ell_2$-strictly singular),
where $p'={p \over {p-1}}$ is the conjugate index to $p$. Better known is that
an operator on $L_p$, $1<p<\infty$, is strictly singular if it is both $\ell_p$-strictly singular and $\ell_2$-strictly singular (and hence $T$ is strictly singular if and only if $T^*$ is strictly singular).  For $p>2$ this is immediate from  \cite{kp}, and Lutz Weis \cite{weis} proved the case $p<2$.

Although Saksman and Tylli did not know a complete
characterization of the weakly compact multiplication operators on
$L(L_p)$, they realized that a classification  must involve
$\ell_p$ and $\ell_2$-strictly singular operators on $L_p$.  This
led Tylli to ask us about possible classifications of the $\ell_p$
and $\ell_2$-strictly singular operators on $L_p$. The $\ell_2$
case is known. It is enough to consider the case $2<p<\infty$.  If
$T$ is an operator on $L_p$, $2<p<\infty$, and $T$ is
$\ell_2$-strictly singular, then it is an easy consequence of the
Kadec-Pe\l czy\'nski dichotomy that $I_{p,2}T$ is compact, where
$I_{p,r}$ is the identity mapping from $L_p$ into $L_r$. But then
by \cite{j}, $T$ factors through $\ell_p$. Tylli then asked
whether the following conjecture is true:

{\bf Tylli Conjecture.} If $T$ is an $\ell_p$-strictly singular
operator on $L_p$, $1<p<\infty$, then   $T$ is in the closure (in
the operator norm) of the operators on $L_p$ that factor through
$\ell_2$. (It is clear that the closure is needed because not all
compact operators on $L_p$, $p\not= 2$, factor through $\ell_2$.)

We then formulated a weaker conjecture:

{\bf Weak Tylli Conjecture.} If $T$ is an $\ell_p$-strictly
singular operator on $L_p$, $1<p<\infty$, and $J:L_p\to
\ell_\infty$ is an isometric embedding, then   $JT$ is in the
closure of the operators from $L_p$ into $\ell_\infty$ that factor
through $\ell_2$.

It is of course evident that an operator on $L_p$, $p\not= 2$,
that satisfies the conclusion of the Weak Tylli Conjecture must be
$\ell_p$-strictly singular. There is a slight subtlety in these
conjectures: while the Tylli Conjecture for $p$ is equivalent to
the Tylli Conjecture for $p'$, it is not at all clear and   is
even false that the Weak Tylli Conjecture for $p$ is equivalent to
the Weak Tylli Conjecture for $p'$.  In fact, we observe in Lemma
\ref{weak} (it is simple) that for $p>2$ the Weak Tylli Conjecture
is true, while the example in Section \ref{example} yields that
the Tylli Conjecture is false for all $p\not= 2$ and the Weak
Tylli Conjecture is false for $p<2$.

There are however some interesting consequences of the Weak Tylli
Conjecture  that are true when $p<2$.  In Theorem
\ref{maintheorem} we prove that if $T$ is an $\ell_p$-strictly
singular operator on $L_p$, $1<p<2$, then $T$ is $\ell_r$-strictly
singular for all $p<r<2$. In view of theorems of Aldous \cite{ald}
(see also \cite{km}) and Rosenthal \cite{ros},  this proves that
if such a $T$ is an isomorphism on a subspace $Z$ of $L_p$, then
$Z$ embeds into $L_r$ for all $r<2$.  The Weak Tylli Conjecture
would imply   that $Z$ is isomorphic to $\ell_2$, but the example
in Section \ref{example} shows that this need not be true.  When
we discovered Theorem \ref{maintheorem}, we thought its  proof
bizarre and assumed that a more straightforward argument would
yield a stronger theorem. The example suggests that in fact the
proof may be ``natural".

In Section  \ref{biased} we discuss   convolution by a biased coin on
$L_p$ of the Cantor group, $1\le p<2$.  We do not know whether such an operator $T$ on $L_p$, $1<p<2$, must satisfy the Tylli Conjecture or the weak Tylli conjecture.  We do prove, however, that if  $T_{|X}$ is an
isomorphism for some reflexive subspace $X$ of $L_p$, $1\le p<2$, then $X$ is
isomorphic to a Hilbert space. This answers an old question of
H.~P.~Rosenthal \cite{rosAltgeld}.

The standard Banach space theory terminology and background we use
can be found in \cite{lt}.

\section{ Weakly compact multiplication operators on $L(L_p)$}

\ \ \ \ We use freely the result \cite[Proposition 2.5]{st2} that
if $A$, $B$ are in $L(X)$ where  $X$ is a  reflexive Banach space
with the approximation property, then the multiplication operator
$L_AR_B$ on $L(X)$ is weakly compact if and only if for every $T$
in $L(X)$, the operator $ATB$ is compact. For completenes, in section \ref{appendix}  we give another proof of this under the weaker assumption that $X$ is reflexive and has the compact approximation property.
This theorem implies that for
such an $X$, $L_AR_B$ is weakly compact on $L(X)$  if and only if
$L_{B^*}R_{A^*}$ is a weakly compact operator on $L(X^*)$.
Consequently, to classify weakly compact multiplication operators
on $L(L_p)$, $1<p<\infty$, it is enough to consider the case
$p>2$. For $p\le r$ we denote the identity operator from $\ell_p$
into $\ell_r$ by $i_{p,r}$. It is immediate from \cite{kp} that an operator $T$  on $L_p$, $2<p<\infty$, is compact if and only if $i_{2,p}$ does not factor through $T$.

\begin{thm}\label{multiplication} Let  $2<p<\infty$ and let $A$, $B$ be bounded linear operators on
$L_p$. Then the multiplication operator $L_AR_B$ on $L(L_p)$ is weakly compact if and only if
one of the following (mutually exclusive) conditions hold.
\begin{itemize}
\item[(a)] $i_{2,p}$ does not factor through $A$ (i.e., $A$ is compact)
\item[(b)] $i_{2,p}$   factors through $A$ but $i_{p,p}$ does not factor through $A$
(i.e., $A$ is $\ell_p$-strictly singular) and $i_{2,2}$ does not factor through $B$
(i.e., $B$ is $\ell_2$-strictly singular)
\item[(c)] $i_{p,p}$ factors through $A$ but $i_{2,p}$ does not factor through $B$
(i.e., $B$ is compact)
\end{itemize}
\end{thm}

\pf The proof is a straightforward application of the Kadec-Pe\l
czy\'nski dichotomy principle \cite{kp}: if $\xn$ is a
semi-normalized (i.e., bounded and bounded away from zero) weakly
null sequence in $L_p$, $2<p<\infty$, then there is a subsequence
which is equivalent to either the unit vector basis of $\ell_p$
or of $\ell_2$ and spans a complemented subspace of $L_p$. Notice
that this immediately implies the ``i.e.'s" in the statement of
the theorem so that (a) and (c) imply that $L_AR_B$ is weakly
compact.  Now assume that (b) holds and let $T$ be in $L(L_P)$.
If $ATB$ is not compact, then there is a normalized weakly null
sequence $\xn$ in $L_p$ so that $ATBx_n$ is bounded away from
zero.  By passing to a subsequence, we may assume that $\xn$ is
equivalent to  either the unit vector basis of $\ell_p$  or of
$\ell_2$.  If $\xn$ is equivalent to  the unit vector basis of
$\ell_p$, then since $TBx_n$ is bounded away from zero, we can
assume by passing to another subsequence that also $TBx_n$ is
equivalent to   the unit vector basis of $\ell_p$ and similarly
for $ATBx_n$, which contradicts the assumption that $A$ is
$\ell_p$-strictly singular. On the other hand, if $\xn$ is
equivalent to  the unit vector basis of $\ell_2$, then since $B$
is $\ell_2$-strictly singular we can assume by passing to a
subsequence that $Bx_n$ is equivalent to   the unit vector basis
of $\ell_p$ and continue as in the previous case to get a
contradiction.

Now suppose that (a), (b), and (c) are all false. If $i_{p,p}$
factors through $A$ and $i_{2,p}$  factors through $B$ then there
is sequence $\xn$ equivalent to the unit vector basis of $\ell_2$ or of  $\ell_p$
so that $Bx_n$ is  equivalent to the unit vector basis of $\ell_2$ or of  $\ell_p$ (of course, only three of the four cases are possible)
and $Bx_n$ spans a complemented subspace of $L_p$.  Moreover,
there is a sequence $\yn$ in $L_p$ so that both $y_n$ and $Ay_n$
are equivalent to the unit vector basis of $\ell_p$. Since $Bx_n$
spans a complemented subspace of $L_p$, the mapping $Bx_n\mapsto
y_n$ extends to a bounded linear operator $T$ on $L_p$ and $ATB$
is not compact.  Finally, suppose that $i_{2,p}$   factors through
$A$ but $i_{p,p}$ does not factor through $A$ and $i_{2,2}$
factors through $B$. Then there is a sequence $\xn$ so that $x_n$
and $Bx_n$ are both equivalent to the unit vector basis of
$\ell_2$ and $Bx_n$  spans a complemented subspace of $L_p$. There
is also a sequence $\yn$  equivalent to the unit vector basis of
$\ell_2$ so that $Ay_n$ is  equivalent to the unit vector basis of $\ell_2$ or of
$\ell_p$.  The mapping $Bx_n\mapsto y_n$ extends to a bounded
linear operator $T$ on $L_p$ and $ATB$ is not compact.
\endpf

It is perhaps worthwhile to restate Theorem \ref{multiplication} in a way that the cases where $L_AR_B$ is weakly compact  are not mutually exclusive.

\begin{thm}\label{multiplicationtwo} Let  $2<p<\infty$ and let $A$, $B$ be bounded linear operators on
$L_p$. Then the multiplication operator $L_AR_B$ on $L(L_p)$ is weakly compact if and only if
one of the following   conditions hold.
\begin{itemize}
\item[(a)]  $A$ is compact
\item[(b)] $A$ is $\ell_p$-strictly singular  and  $B$ is $\ell_2$-strictly singular
\item[(c)]  $B$ is compact
\end{itemize}
\end{thm}

\section{ $\ell_p$-strictly singular operators on $L_p$}

\ \ \ \ We recall the well known

\begin{lm}\label{ellpsequence} Let $W$ be a bounded convex symmetric subset of $L_p$, $1\le
p\not=2<\infty$. The following are equivalent:
\begin{enumerate}
\item[1.] No sequence  in $W$   equivalent to the unit vector basis for $\ell_p$ spans
a complemented subspace of $L_p$.
\item[2.] For every $C$ there exists $n$ so that no length $n$ sequence in $W$ is
$C$-equivalent to the unit vector basis of $\ell_p^n$.
\item[3.] For each $\ep>0$ there is $M_\ep$ so that
$W\subset \ep B_{L_p} + M_\ep B_{L_\infty}$.
\item[4.] $|W|^p$ is  uniformly integrable; i.e.,
$\lim_{t\downarrow 0} \sup_{x\in W} \sup_{\mu(E)< t} \|\1_E x\|_p
=0$.
\end{enumerate}
\end{lm}

When $p=1$, the assumptions that $W$ is convex and and $W$
symmetric are not needed, and the conditions in Lemma
\ref{ellpsequence} are equivalent to the non weak compactness of
the weak closure of $ {W}$. This case is essentially proved in
\cite{kp} and proofs  can also be found in books; see, e.g.,
\cite[Theorem 3.C.12]{woj}). (Condition (3) does not appear in
\cite{woj}, but it is easy to check the equivalence of (3) and
(4).  Also, in the proof in \cite[Theorem 3.C.12]{woj}) that not
(4) implies not (1), Wojtaszczyk only constructs a basic sequence
in $W$ that is equivalent to the unit vector basis for $\ell_1$;
however, it is clear that the constructed basic sequence spans a
complemented subspace.)

 For $p>2$, Lemma
\ref{ellpsequence} and stronger versions of condition (1) can be
deduced from \cite{kp}. For $1<p<2$, one needs to modify slightly
the proof in \cite{woj} for the case $p=1$.  The only essential
modification comes in the proof that not (4) implies not (1), and
this is where it is needed that $W$ is convex and symmetric.  Just
as in \cite{woj}, one shows that not (4) implies that there is a
sequence $\xn$ in $W$ and a sequence $\En$ of disjoint measurable
sets so that $\inf \|1_{E_n} x_n\|_p >0$.  By passing to a
subsequence, we can assume that $\xn$ converges weakly to, say,
$x$.  Suppose first that $x=0$. Then by passing to a further
subsequence, we may assume that $\xn$ is a small perturbation of a
block basis of the Haar basis for $L_p$ and hence is an
unconditionally basic sequence. Since $L_p$ has type $p$, this
implies that there is a constant $C$ so that for all sequences
$\an$ of scalars, $\|\sum a_n x_n\|_p\le C(\sum |a_n|^p)^{1/p}$.
Let $P$ be the norm one projection from $L_p$ onto the closed
linear span $Y$ of the disjoint sequence $ \Enxn$. Then $Px_n$ is
weakly null in a space isometric to $\ell_p$ and $\|Px_n\|_p$ is
bounded away from zero, so there is a subsequence
$\{Px_{n(k)}\}_{k=1}^\infty$ which is equivalent to the unit
vector basis for $\ell_p$ and whose closed span  is the range of a
projection $Q$ from $Y$.  The projection $QP$ from $L_p$ onto the
the closed span of $\{Px_{n(k)}\}_{k=1}^\infty$ maps $x_{n(k)}$ to
$Px_{n(k)}$ and, because of the upper $p$ estimate on $\xnk$, maps
the closed span of $\xnk$ isomorphically onto the closed span of
   $\{Px_{n(k)}\}_{k=1}^\infty$.  This yields that $\xnk$ is
equivalent to the unit vector basis for $\ell_p$ and spans a
complemented subspace.  Suppose now that the weak limit $x$ of
$\xn$ is not zero. Choose a subsequence $\xnk$ so that $\inf
\|1_{E_{n(2k+1)}}(x_{n(2k)}-x_{n(2k+1)})\|_p>0$ and replace $\xn$
with $\{{{x_{n(2k)}-x_{n(2k+1)}}\over 2}\}_{k=1}^\infty$ in the
argument above.

Notice that the argument outlined above gives that if $\xn$ is a
sequence in $L_p$, $1< p\not=2<\infty$, which is equivalent to the
unit vector basis of $\ell_p$, then there is a subsequence $\yn$
 whose closed linear span in $L_p$ is complemented.  This
is how one proves that the identity on $\ell_p$ factors through
any operator on $L_p$ which is not $\ell_p$-strictly singular.

The Weak Tylli Conjecture for $p>2$ is an easy consequence of the
following lemma.

\begin{lm}\label{weak} Let $T$ be an  operator from a
$\mathcal{L}_{1}$ space $V$ into $L_p$, $1<p<2$,  so that $W:=TB_{V}$ satisfies
condition (1) in Lemma \ref{ellpsequence}. Then for each
$\ep>0$  there is an operator $S:V\to L_2$ so that   $\|T-I_{2,p}S\| <\ep$.
 \end{lm}
\pf Let $\ep>0$.  By condition (3) in Lemma \ref{ellpsequence},
for each norm one vector $x$ in $V$ there is a vector $Ux$ in
$L_2$ with $\|Ux\|_2\le \|Ux\|_\infty\le  M_\ep$ and
$\|Tx-Ux\|_p\le \ep$.  By the definition of $\mathcal{L}_{1}$
space, we can write $V$ as a directed union $\cup_\alpha E_\alpha$
of finite dimensional spaces that are uniformly isomorphic to
$\ell_1^{n_\alpha}$, $n_\alpha=\dim E_\alpha$,   and let
$(x^\alpha_i)_{i=1}^{n_\alpha}$ be norm one vectors in $E_\alpha$
which are, say, $\lambda$-equivalent to the unit vector basis for
$\ell_1^{n_\alpha}$ with $\lambda$ independent of $\alpha$.  Let
$U_\alpha$ be the linear extension to $E_\alpha$ of the mapping
$x^\alpha_i\mapsto Ux^\alpha_i$, considered as an operator into
$L_2$. Then $\|T_{|E_\alpha}- I_{2,p}U_\alpha\|\le \lambda \ep$
and $\|U_\alpha\|\le \lambda M_\ep$.  A standard Lindenstrauss
compactness argument produces an operator $S:V\to L_2$ so that
$\|S\|\le \lambda M_\ep$ and $\|T-I_{2,p}S\|\le \lambda \ep$.
Indeed, extend $U_\alpha$ to all of $V$ by letting $U_\alpha x=0 $
if $x\not\in E_\alpha$. The net $T_\alpha$ has a subnet $S_\beta$
so that for each $x$ in $V$, $S_\beta x$ converges weakly in
$L_2$; call the limit $Sx$. It is easy to check that $S$ has the
properties claimed.
\endpf

\begin{thm}\label{weaktyllithm} Let $T$ be an $\ell_p$-strictly singular operator on
$L_p$, $2<p<\infty$, and let $J$ be an isometric embedding of $L_p$ into an injective $Z$.
  Then for each $\ep>0$ there is an operator $S:L_p\to Z$ so that $S$ factors
through $\ell_2$ and $\|JT-S\| <\ep$.
\end{thm}
\pf Lemma \ref{weak} gives the conclusion when $J$ is the adjoint
of a quotient mapping from $\ell_1$ or $L_1$ onto $L_{p'}$. The
general case then follows from the injectivity of $Z$.
\endpf

The next proposition, when souped-up via ``abstract nonsense" and
known results,
 gives our main
result about
$\ell_p$-strictly singular operators on $L_p$.  Note that it shows that an
$\ell_p$-strictly singular operator on $L_p$, $1<p<2$, cannot be the identity on
the span
of a sequence of
$r$-stable independent random variables for any $p<r<2$.  We do not know another way of
proving even this special case of our main result.

\begin{pr}\label{proposition1} Let $T$ be an $\ell_p$-strictly singular operator on
$L_p$, $1<p<2$. If $X$ is a subspace of $L_p$ and $T_{|X}=aI_X$ with $a\not= 0$, then
$X$ embeds into $L_s$ for all $s<2$.
\end{pr}
\pf By making  a change of density, we can by \cite{jj} assume
that $T$ is also a bounded linear operator on $L_2$, so assume,
without loss of generality, that $\|T\|_p\vee \|T\|_2 = 1$, so
that, in particular, $a\le 1$.  Lemma \ref{ellpsequence} gives for
each $\epsilon>0$ a constant $M_\epsilon$ so that
\begin{equation}\label{eqprop1}
T B_{L_p}\subset \epsilon B_{L_p} +  M_\epsilon B_{L_2}.
\end{equation}

Indeed, otherwise condition (1) in Lemma \ref{ellpsequence} gives
a bounded sequence $\xn$ in $L_p$ so that $\Txn$   is equivalent
to the unit vector basis of $\ell_p$.  By passing to a subsequence
of differences of $\xn$, we can assume, without loss of
generality, that $\xn$ is a small perturbation of a block basis of
the Haar basis for $L_p$ and hence is an unconditionally basic
sequence. Since $L_p$ has type $p$, the sequence $\xn$ has an
upper $p$ estimate, which means that there is a constant $C$ so
that for all sequences $\an$ of scalars, $\| \sum a_n x_n\| \le
C\| (\sum |a_n|^p)^{1/p}\|$. Since $\Txn$ is equivalent to the
unit vector basis of $\ell_p$, $\xn$ also has a lower $p$ estimate
and hence $\xn$ is equivalent to the unit vector basis of
$\ell_p$. This contradicts the $\ell_p$ strict singularity of $T$.

Iterating this we get for every $n$ and $0<\epsilon<1/2$
\[
a^n B_X \subset T^n B_{L_p}\subset \epsilon^n B_{L_p} +  2M_\epsilon B_{L_2}
\]
or, setting $A:=1/a$,
\[
B_X \subset  A^n\epsilon^n B_{L_p} +  2 A^n M_\epsilon B_{L_2}.
\]

For $f$ a unit vector in $X$ write $f=f_n+g_n$ with $\|f_n\|_2\le  2 A^n M_\epsilon$
and $\|g_n\|_p\le (A\epsilon)^n$. Then $f_{n+1}-f_n=g_n-g_{n+1}$, and since
evidently $f_n$ can be chosen to be of the form $(f\vee -k_n)\wedge k_n$ (with
appropriate interpretation when the set $[f_n=\pm k_n]$ has  positive measure), the
choice of
$f_n$, $g_n$ can be made so that
\[
\|f_{n+1}-f_n\|_2\le \|f_{n+1}\|_2\le 2M_\epsilon A^{n+1}
\]
\[
\|g_n-g_{n+1}\|_p\le \|g_n\|_p\le (A\epsilon)^n.
\]

For $p<s<2$ write ${1\over s}={\theta\over 2} +{{1-\theta}\over p}$.  Then
\[
\|f_{n+1}-f_n\|_s  \le \|f_{n+1}-f_n\|_2^\theta \|g_n-g_{n+1}\|_p^{1-\theta}
\le(2M_\epsilon A)^\theta (A\epsilon^{1-\theta})^n
\]
which is summable if $\epsilon^{1-\theta}<1/A$.   But $\|f-f_n\|_p\to 0$ so
$f=f_1 + \sum_{n=1}^\infty f_{n+1} - f_n$ in $L_p$ and hence also in $L_s$ if
$\epsilon^{1-\theta}<1/A$. So for some constant $C_s$ we get for all $f\in X$ that
$\|f\|_p\le \|f\|_s\le C_s\|f\|_p$. \qed

\medskip
We can now prove our main theorem. For background on ultrapowers
of Banach spaces, see \cite[Chapter 8]{djt}.

\begin{thm}\label{maintheorem} Let $T$ be an $\ell_p$-strictly singular operator on
$L_p$, $1<p<2$. If $X$ is a subspace of $L_p$ and $T_{|X}$ is an isomorphism, then $X$ embeds
into $L_r$ for all $r<2$.
\end{thm}
\pf In view of Rosenthal's theorem \cite{ros}, it is enough to
prove that $X$ has type $s$ for all $s<2$.  By virtue of of the
Krivine-Maurey-Pisier theorem, \cite{kr} and \cite{mp} (or,
alternatively, Aldous' theorem, \cite{ald} or \cite{km}), we only
need to check that for $p<s<2$, $X$ does not contain almost
isometric copies of $\ell_s^n$ for all $n$. (To apply the
Krivine-Maurey-Pisier theorem we use that the second condition in
Lemma \ref{ellpsequence}, applied to the unit ball of $X$, yields
that $X$ has type $s$ for some $p<s\le 2$). So suppose that for
some $p<s<2$, $X$   contains almost isometric copies of $\ell_s^n$
for all $n$. By applying Krivine's theorem \cite{kr} we get for
each $n$ a sequence $(f^n_i)_{i=1}^n$ of unit vectors in $X$ which
is $1+\epsilon$-equivalent to the unit vector basis for $\ell_s^n$
and, for some constant $C$ (which we can take independently of
$n$), the sequence $(CTf^n_i)_{i=1}^n$ is also
$1+\epsilon$-equivalent to the unit vector basis for $\ell_s^n$.
By replacing $T$ by $CT$, we might as well assume that $C=1$.  Now
consider an ultrapower $T_\calU$, where $\calU$ is a free
ultrafilter on the natural numbers.  The domain and codomain of
$T_\calU$ is the (abstract) $L_p$ space $(L_p)_\calU$, and
$T_\calU$ is defined by $T_\calU(f_1,f_2,\dots) = (Tf_1,
Tf_2,\dots)$ for any (equivalence class of a) bounded sequence
$(f_1,f_2,\dots)$.  It is evident that $T_\calU$ is an isometry on
the ultraproduct of span $(f^n_i)_{i=1}^n$; $n=1,2,\dots$,  and
hence $T_\calU$ is an isometry on a subspace of $(L_p)_\calU$
which is isometric to $\ell_s$. Since condition 2 in Lemma
\ref{multiplication} is obviously preserved when taking an
ultrapower of a set, we see that $T_\calU$ is $\ell_p$-strictly
singular.  Finally, by restricting $T_\calU$ to a suitable
subspace, we get an $\ell_p$-strictly singular operator $S$ on
$L_p$ and a subspace $Y$ of $L_p$ so that $Y$ is isometric to
$\ell_s$ and $S_{|Y}$ is an isometry.  By restricting the domain
of $S$, we can assume that $Y$ has full support and the functions
in $Y$ generate the Borel sets.  It then follows from the
Plotkin-Rudin theorem \cite{pl}, \cite{ru} (see \cite[Theorem
1]{kk}) that $S_{|Y}$ extends to an isometry $W$ from $L_p$ into
$L_p$. Since any isometric copy of $L_p$ in $L_p$ is norm one
complemented (see \cite[\S 17]{la}), there is a norm one operator
$V:L_p\to L_p$ so that $VW=I_{L_p}$. Then $VS_{|Y}=I_Y$ and $VS$
is $\ell_p$-strictly singular, which contradicts Proposition
\ref{proposition1}.
\endpf

\begin{re} The   $\ell_1$-strictly singular operators on
$L_1$ also form an interesting class.  They are the weakly compact
operators on $L_1$.  In terms of factorization, they are just the
closure in the operator norm of the integral operators on $L_1$
(see, e.g., the proof of Lemma \ref{weak}).
\end{re}

\section{The example}\label{example}

\ \ \ \ Rosenthal \cite{rosxp} proved that if $\xn$ is a sequence
of three valued, symmetric, independent random variables, then for
all $1<p<\infty$, the closed span in $L_p$ of $\xn$ is
complemented by means of the orthogonal projection $P$, and
$\|P\|_p$ depends only on $p$, not on the specific sequence $\xn$.
Moreover, he showed that if $p>2$, then for any sequence $\xn$ of
symmetric, independent random variables in $L_p$, $\|\sum x_n\|_p$
is equivalent (with constant depending only on $p$) to $(\sum
\|x_n\|_p^p)^{1/p}\vee (\sum \|x_n\|_2^2)^{1/2}$. Thus if $\xn$ is
normalized in $L_p$, $p>2$, and  $w_n:=\|x_n\|_2$, then $\|\sum
a_n x_n\|_p$ is equivalent to $\|\an\|_{p,w}:=(\sum
|a_n|^p)^{1/p}\vee (\sum |a_n|^2w_n^2)^{1/2}$. The completion of
the finitely non zero sequences of scalars under the norm
$\|\cdot\|_{p,w}$ is called $X_{p,w}$. It follows that if $w=\wn$
is any sequence of numbers in $[0,1]$, then $X_{p,w}$ is
isomorphic to a complemented subspace of $L_p$. Suppose now that
$w=\wn$ and $v=\vn$ are two such sequences of weights and $v_n\ge
w_n$, then the diagonal operator $D$ from $X_{p,w}$ to $X_{p,v}$
that sends the $n$th unit vector basis vector $e_n$ to
${{w_n}\over{v_n}} e_n$ is contractive and it is more or less
obvious that $D$  is $\ell_p$-strictly singular if
${{w_n}\over{v_n}}\to 0$ as $n\to \infty$. Since $X_{p,w}$ and
$X_{p,v}$ are isomorphic to complemented subspaces of $L_p$, the
adjoint operator $D^*$ is $\ell_{p'}$ strictly singular and
(identifying $X_{p,w}^*$ and $X_{p,v}^*$ with subspaces of
$L_{p'}$) extends to a $\ell_{p'}$ strictly singular operator on
$L_{p'}$.  Our goal in this section is produce weights $w$ and $v$
so that $D^*$ is an isomorphism on a subspace of $X_{p,v}^*$ which
is not isomorphic to a Hilbert space.

For all $0<r<2$ there is a positive constant $c_r$ such that
\[
|t|^r=c_r\int_0^\infty \frac{1-\cos tx}{x^{r+1}} dx
\]
for all $t\in\R$. It follows that for any closed interval
$[a,b]\subset (0,\infty)$ and for all $\e>0$ there are
$0<x_1<x_2<\dots<x_{n+1}$ such that $\max_{1\le j\le
n}\Big|\frac{x_{j+1}-x_j}{x_j^{r+1}}\Big|\le\e$ and
\begin{equation}\label{eq:appr}
\Big|c_r\sum_{j=1}^n\frac{x_{j+1}-x_j}{x_j^{r+1}}(1-\cos
tx_j)-|t|^r\Big|<\e
\end{equation}
for all $t$ with $|t|\in[a,b]$.

Let $0<q<r<2$ and define $v_j$ and $a_j$, $j=1,\dots,n$, by
\[
v_j^{\frac{2q}{2-q}}=c_r\frac{x_{j+1}-x_j}{x_j^{r+1}} \ , \ \ \ \
\ \frac {a_j}{v_j^{\frac{2}{2-q}}}=x_j.
\]
 Let $Y_j$, $j=1,\dots,n$, be independent, symmetric, three valued
 random variables such that $|Y_j|=v_j^{\frac{-2}{2-q}}\1_{B_j}$
 with ${\rm Prob}(B_j)=v_j^{\frac{2q}{2-q}}$, so that in particular
 $\|Y_j\|_q=1$ and $v_j=\|Y_j\|_q/\|Y_j\|_2$. Then the characteristic function of $Y_j$ is
 \[
 \varphi_{Y_j}(t)=1-v_j^{\frac{2q}{2-q}}+v_j^{\frac{2q}{2-q}}\cos
 (tv_j^{\frac{-2}{2-q}})=1-v_j^{\frac{2q}{2-q}}(1-\cos
 (tv_j^{\frac{-2}{2-q}}))
 \]
 and
\begin{equation}\label{eq:phi}
\begin{array}{rl}
\varphi_{\sum a_jY_j}(t)= &
\prod_{j=1}^n(1-v_j^{\frac{2q}{2-q}}(1-\cos
 (ta_jv_j^{\frac{-2}{2-q}})))\\
 =&
\prod_{j=1}^n(1-c_r\frac{x_{j+1}-x_j}{x_j^{r+1}}(1-\cos
 (tx_j)))
\end{array}
\end{equation}
 To evaluate this product we use the estimates
 on $\frac{x_{j+1}-x_j}{x_j^{r+1}}$ to deduce that, for each $j$
 \[
 \begin{array}{rl}
 |\log(1-c_r\frac{x_{j+1}-x_j}{x_j^{r+1}}(1-\cos
 (tx_j)))+&c_r\frac{x_{j+1}-x_j}{x_j^{r+1}}(1-\cos
 (tx_j))|\\
 &\le C\e c_r^2\frac{x_{j+1}-x_j}{x_j^{r+1}}(1-\cos
 (tx_j))
 \end{array}
 \]
 for some absolute $C<\infty$. Then, by (\ref{eq:appr}),
 \[
 \begin{array}{rl}
 |\sum_{j=1}^n\log(1-c_r\frac{x_{j+1}-x_j}{x_j^{r+1}}(1-\cos
 (tx_j)))+&c_r\sum_{j=1}^n\frac{x_{j+1}-x_j}{x_j^{r+1}}(1-\cos
 (tx_j))|\\
 &\le C\e c_r(\e+b^r).
 \end{array}
 \]
Using (\ref{eq:appr}) again we get
\[
|\sum_{j=1}^n\log(1-c_r\frac{x_{j+1}-x_j}{x_j^{r+1}}(1-\cos
 (tx_j)))+|t|^r|\le (C+1)\e(\e+b^r)
 \]
 (assuming as we may that $b\ge 1$), and from (\ref{eq:phi}) we get
 \[
 \varphi_{\sum a_jY_j}(t)=(1+O(\e))\exp(-|t|^r)
 \]
for all $|t|\in [a,b]$, where the function hiding under the $O$
notation depends on $r$ and $b$ but on nothing else. It follows
that, given any $\eta>0$, one can find $a,b$ and $\e$, such that
for the corresponding $\{a_j,Y_j\}$ there is a symmetric
$r$-stable $Y$ (with characteristic function $e^{-|t|^r}$)
satisfying
\[
\|Y-\sum_{j=1}^n a_j Y_j\|_q\le\eta.
\]
This follows from classical translation of various convergence
notions; see e.g.   \cite[p. 154]{rosII}.

 Let now $0<\delta<1$. Put
$w_j=\delta v_j$, $j=1,\dots,n$, and let $Z_j$, $j=1,\dots,n$, be
independent, symmetric, three valued
 random variables such that $|Z_j|=w_j^{\frac{-2}{2-q}}\1_{C_j}$
 with ${\rm Prob}(C_j)=w_j^{\frac{2q}{2-q}}$, so that in particular
 $\|Z_j\|_q=1$ and $w_j=\|Z_j\|_q/\|Z_j\|_2$. In a similar manner
 to the argument above we get that,

 \begin{equation*}
\begin{array}{rl}
\varphi_{\sum \delta a_jZ_j}(t)= &
\prod_{j=1}^n(1-w_j^{\frac{2q}{2-q}}(1-\cos
 (t\delta a_j w_j^{\frac{-2}{2-q}})))\\
 = &
 \prod_{j=1}^n(1-\delta^{\frac{2q}{2-q}}v_j^{\frac{2q}{2-q}}(1-\cos
 (t\delta^{\frac{-q}{2-q}} a_j v_j^{\frac{-2}{2-q}})))\\
 =
 & (1+O(\e))\exp(-\delta^{\frac{q(2-r)}{2-q}}|t|^r)
\end{array}
\end{equation*}
 for all $|t|\in
[\delta^{\frac{q}{2-q}}a,\delta^{\frac{q}{2-q}}b]$, where the $O$
now depends also on $\delta$.

Assuming $\delta^{\frac{q(2-r)}{2-q}}>1/2$ and for a choice of
$a,b$ and $\e$, depending on $\delta, r, q$ and $\eta$ we get that
there is a symmetric $r$-stable random variable $Z$ (with
characteristic function $e^{-\delta^\frac{q(2-r)}{2-q}|t|^r}$)
such that
\[
\|Z-\sum_{j=1}^n \delta a_j Z_j\|_q\le\eta.
\]
Note that the ratio between the $L_q$ norms of $Y$ and $Z$ are
bounded away from zero and infinity by universal constants and
each of these norms is also universally bounded away from zero.
Consequently, if $\e$ is small enough the ratio between the $L_q$
norms of $\sum_{j=1}^n a_j Y_j$ and $\sum_{j=1}^n \delta a_j Z_j$
are bounded away from zero and infinity by universal constants.

Let now $\delta_i$ be any sequence decreasing to zero and $r_i$
any sequence such that $q<r_i\uparrow 2$ and satisfying
$\delta_i^{\frac{q(2-r_i)}{2-q}}>1/2$. Then, for any sequence
$\e_i\downarrow 0$ we can find two sequences of symmetric,
independent, three valued random variables $\{Y_i\}$ and
$\{W_i\}$, all normalized in $L_q$, with the following additional
properties:

\begin{itemize}
\item put $v_j=\|Y_j\|_q/\|Y_j\|_2$ and $w_j=\|Z_j\|_q/\|Z_j\|_2$.
Then there are disjoint finite subsets of the integers $\sigma_i$,
$i=1,2,\dots$, such that $w_j=\delta_i v_j$ for $j\in \sigma_i$.
\item There are independent random variables $\{\bar Y_i\}$ and $\{\bar Z_i\}$
with $\bar Y_i$ and $\bar Z_i$ $r_i$ stable with bounded, from
zero and infinity, ratio of $L_q$ norms and there are coefficients
$\{a_j\}$ such that
\[
\|\bar Y_i-\sum_{j\in\sigma_i}a_jY_j\|_q<\e_i \ \ \ \mbox {and}\ \
\ \|\bar Z_i-\sum_{j\in\sigma_i}\delta_ia_jZ_j\|_q<\e_i.
\]
\end{itemize}

From \cite{rosxp} we know that the spans of $\{Y_j\}$ and
$\{Z_j\}$ are complemented in $L_q$, $1<q<2$, and the dual spaces
are naturally isomorphic to the spaces $X_{p,\{v_j\}}$ and
$X_{p,\{w_j\}}$ respectively; both the isomorphism constants and
the complementation constants depend only on $q$. Here $p=q/(q-1)$
and
\[
\|\{\alpha_j\}\|_{X_{p,\{u_j\}}}=\max\{(\sum|\alpha_j|^p)^{1/p},(\sum
u_j^2\alpha_j^2)^{1/2}\}.
\]
Under this duality the adjoint $D^*$ to the operator $D$ that
sends $Y_j$ to $\delta_i Z_j$ for $j\in\sigma_i$ is formally the
same diagonal operator between $X_{p,\{w_i\}}$ and
$X_{p,\{v_i\}}$. The relation $w_j=\delta_i v_j$ for $j\in
\sigma_i$ easily implies that this is a bounded operator.
$\delta_i\to 0$ implies that this operator is $\ell_q$ strictly
singular.
If $\e_i\to 0$ fast enough, $D^*$ preserves a copy of ${\rm
span}\{\bar Y_i\}$. Finally, if $r_i$ tend to 2 not too fast this
span is not isomorphic to a Hilbert space. Indeed, let $1\le
s_j\uparrow 2$ be arbitrary and let $\{n_j\}_{j=1}^\infty$ be a
sequence of positive integers with $n_j^{\frac{1}{s_j}-\frac12}\ge
j$, $j=1,2,\dots$, say. For $1\le k\le n_j$, put
$r_{n_1+n_2+\dots+n_{j-1}+k}=s_j$. Then the span of
$\{Y_i\}_{i=n_1+\dots+n_{j-1}+1}^{n_1+\dots+n_j}$ is isomorphic,
with constant independent of $j$, to $\ell_{s_j}^{n_j}$ and this
last space is of distance at least $j$ from a Euclidean space.

It follows that if $J:L_q\to \ell_\infty$ is an isometric
embedding, then $JD^*$
 cannot be arbitrarily
approximated by an operator which factors through a Hilbert space,
and hence the Weak Tylli Conjecture is false in the range $1<q<2$.

\section{Convolution by a biased coin}\label{biased}

\ \ \ \ In this section we regard $L_p$ as $L_p(\Delta)$, where
$\Delta=\{-1,1\}^\N$ is the Cantor group and the measure is the
Haar measure $\mu$ on $\Delta$; i.e., $\mu=\prod_{n=1}^\infty
\mu_n$, where $\mu_n(-1)=\mu_n(1)= 1/2$.  For $0<\e<1$, let
$\nu_\e$ be the $\e$- biased coin tossing measure; i.e.,
$\nu_\e=\prod_{n=1}^\infty \nu_{\e,n}$, where
$\nu_{\e,n}(1)=\frac{1+\e}{2}$ and
$\nu_{\e,n}(-1)=\frac{1-\e}{2}$.  Let $T_\e$ be convolution by
$\nu_\e$, so that for a $\mu$-integrable function $f$ on $\Delta$,
$(T_\e f)(x)=(f*\nu_\e)(x)=\int_\Delta f(xy)\, d\nu_\e (y)$.  The
operator $T_\e$ is a contraction on $L_p$ for all $1\le p\le
\infty$. Let us recall how $T_\e$ acts on the characters on
$\Delta$.  For $t=\tn \in \Delta$, let $r_n(t)=t_n$.  The
characters on $\Delta$ are finite products of these {\it
Rademacher functions} $r_n$ (where the void product is the
constant one function).  For $A$ a finite subset of $\N$, set
$w_A=\prod_{n\in A} r_n$ and let $W_n$ be the linear span of
$\{w_A: |A|=n\}$. Then $T_\e w_A=\e^{|A|}w_A$.

 We are interested in studying $T_\e$ on $L_p$, $1\le p
<2$. The background we mention below is all contained in Bonami's
paper \cite{bon} (or see \cite{rosAltgeld}). On $L_p$, $1<p<2$,
$T_\e$ is $\ell_p$-strictly singular; in fact, $T_\e$ even maps
$L_p$ into $L_r$ for some $r=r(p,\e)>p$. Indeed, by interpolation
it is sufficient to check that $T_\e $ maps $L_s$ into $L_2$ for
some $s=s(\e)<2$. But there is a constant $C_s$ which tends to $1$
as $s\uparrow 2$ so that for all $f\in W_n$, $\|f\|_2\le C_s^n
\|f\|_s$ and the orthogonal projection $P_n$ onto (the closure of)
$W_n$ satisfies $\|P_n\|_p\le C_s^n$.  From this it is easy to
check that  if $\e C_s^2<1$,   then $T_\e$ maps $L_s$ into $L_2$.
We remark in passing that Bonami \cite{bon} found for each $p$ (including $p\ge 2$) and $\e$ the largest value of $r=r(p,\e)$ such that $T_\e$ maps $L_p$ into $L_r$.

Thus Theorem \ref{maintheorem} yields that if $X$ is a subspace of
$L_p$, $1<p<2$, and $T_\e$ (considered as an operator from $L_p$
to $L_p$) is an isomorphism on $X$, then $X$ embeds into $L_s$ for
all $s<2$. Since, as we mentioned above, $T_\e$ maps $L_s$ into
$L_2$ for some $s<2$, it then follows from an argument in
\cite{rosAltgeld} that $X$ must be isomorphic to a Hilbert space. (Actually, as we show after the proof, Lemma \ref{rosLemma} is strong enough that we can prove Theorem \ref{rosProblem} without using Theorem \ref{maintheorem}.)
Since \cite{rosAltgeld} is not generally available, we repeat
Rosenthal's argument in Lemma \ref{rosLemma} below.

Now $T_\e$ is not $\ell_1$ strictly singular on $L_1$.
Nevertheless, we still get that if $X$ is a reflexive subspace of
$L_1$ and $T_\e$ (considered as an operator from $L_1$ to $L_1$)
is an isomorphism on $X$, then $X$ is isomorphic to a Hilbert
space. Indeed, Rosenthal showed (see Lemma \ref{rosLemma}) that
then there is another subspace $X_0$ of $L_1$ which is isomorphic
to $X$ so that $X_0$ is contained in $L_p$ for some $1<p<2$, the
$L_p$ and $L_1$ norms are equivalent on $X_0$, and $T_\e$ is an
isomorphism on $X_0$.  This implies that as an operator on $L_p$,
$T_\e$ is an isomorphism on $X_0$ and hence  $X_0$ is isomorphic
to a Hilbert space. (To apply Lemma \ref{rosLemma}, use the fact
\cite{ros} that if $X$ is a relexive subspace of $L_1$, then $X$
embeds into $L_{p}$ for some $1<p<2$.)

We summarize this discussion in the first sentence of Theorem \ref{rosProblem}.
The case $p=1$ solves Problem B from Rosenthal's 1976 paper
\cite{rosAltgeld}.

\begin{thm}\label{rosProblem} Let $1\le p < 2$, let $0<\e<1$, and
let $T_\e$ be considered as an operator on $L_p$.  If $X$ is a
reflexive subspace of $L_p$ and the restriction of $T_\e$ to $X$
is an isomorphism, then $X$ is isomorphic to a Hilbert space. Moreover, if $p>1$, then $X$ is complemented in $L_p$.
\end{thm}

We now prove Rosenthal's lemma \cite[proof of Theorem 5]{rosAltgeld} and defer the proof of the ``moreover" statement in
Theorem \ref{rosProblem} until after the proof of the lemma.
.

\begin{lm}\label{rosLemma} Suppose that $T$ is an operator on
$L_p$, $1\le p<r<s <2$, $X$ is a subspace of $L_p$ which is
isomorphic to a subspace of $L_s$, and $T_{|X}$ is an isomorphism.
Then there is another subspace $X_0$ of $L_p$ which is isomorphic
to $X$ so that $X_0$ is contained in $L_r$, the $L_r$ and $L_p$
norms are equivalent on $X_0$, and $T$ is an isomorphism on $X_0$.
\end{lm}

\pf We want to find a measurable set $E$ so that
\begin{itemize}
\item[(1)] $X_0:= \{\1_E x : x\in X\}$ is isomorphic to $X$,
\item[(2)] $X_0\subset L_r$,
\item[(3)] $T_{|X_0}$ is an isomorphism.
\end{itemize}
(We did not say that   $\|\cdot\|_p$ and $\|\cdot\|_r$, are
equivalent on $X_0$ since that follows formally from the closed
graph theorem. The isomorphism $X\to X_0$ guaranteed by (a) is of
course the mapping $x\mapsto \1_E x$.)

Assume, without loss of generality, that $\|T\|=1$. Take $a>0$ so
that $\|Tx\|_p \ge {a} \| x\|_p$ for all $x$ in $X$. Since
$\ell_p$ does not embed into $L_s$ we get from (4) in Lemma
\ref{ellpsequence} that there is $\eta>0$ so that if $E$ has
measure larger than $1-\eta$, then $\|\1_{\sim E}x\|_p \le
{\frac{a}{2}}\|x\|_p$ for all $x$ in $x$. Obviously (1) and (3)
are satisfied for any such $E$.  It is proved in \cite{ros} that
    there is strictly positive   $g$ with $\|g\|_1=1$ so that
$\frac {x}{g}$ is in $L_r$ for all $x$ in $X$.  Now simply choose
$t<\infty$ so that $E:=[g<t]$ has measure at least $1-\eta$; then
$E$ satisfies (1), (2), and (3).
\endpf

Next we remark how to avoid using Theorem \ref{maintheorem} in proving Theorem \ref{rosProblem}.  Suppose that  $T_\e$  is an isomorphism on a reflexive subspace $X$ of $ L_p$, $1\le p<2$.  Let $s$ be the supremum of those $r\le 2$ such that $X$ is isomorphic to a subspace of $L_r$, so $1<s\le 2$.  It is sufficient to show that $s=2$. But if $s<2$, we get from the interpolation formula that if $r<s$ is sufficiently close to $s$, then $T_\e$ maps $L_r$ into $L_t$ for some $t>s$ and hence, by Lemma  \ref{rosLemma}, $X$ embeds into $L_t$.

Finally we prove the ``moreover" statement in Theorem \ref{rosProblem}. We now know that $X$ is isomorphic to a Hilbert space.  In the proof of Lemma \ref{rosLemma}, instead of using Rosenthal's result from  \cite{ros}, use Grothendieck's theorem \cite[Theorem 3.5]{djt}, which implies that  there is strictly positive   $g$ with $\|g\|_1=1$ so that
$\frac {x}{g}$ is in $L_2$ for all $x$ in $X$.  Choosing $E$ the same way as in the proof of Lemma \ref{rosLemma} with $T:=T_\e$, we get that (1), (2), and (3) are true with $r=2$. Now the $L_2$ and $L_p$ norms are equivalent on both $X_0$ and on $T_\e X_0$. But it is clear that the only way that $T_\e$ can be an isomorphism on a subspace $X_0$ of $L_2$ is for the orthogonal projection $P_n$ onto the closed span of $W_k$, $0\le k\le n$, to be an isomorphism on $X_0$ for some finite $n$. But then also in the $L_p$ norm the restriction of $P_n$ to $X_0$ is an isomorphism because the $L_p$ norm and the $L_2$ norm are equivalent on the span of $W_k$, $0\le k\le n$ and $P_n$ is bounded on $L_p$ (since $p>1$).  It follows that the operator $S:= P_n \circ \1_E$ on $L_p$ maps $X_0$ isomorphically onto a complemented subspace of $L_p$, which implies that $X_0$ is also complemented in $L_p$.

We conclude this section with the open problem that started us thinking about $\ell_p$-strictly singular operators.

\begin{problem}\label{problem} Let $1<p<2$ and $0<\e<1$.  On $L_p(\Delta)$, does
$T_\e$ satisfy the conclusion of the Tylli Conjecture or the Weak
Tylli Conjecture?
\end{problem}

Of course, the answer to Problem \ref{problem} is ``yes" when
$\e=\e(p)$ is sufficiently small, since then $T_\e$ maps $L_p$
into $L_2$.

\section{Appendix}\label{appendix}

\ \ \  In this appendix we prove a theorem that is essentially due to Saksman and Tylli.  The only novelty is that we assume the compact approximation property rather than the approximation property.

\begin{thm}\label{SaksTylli} Let $X$ be a reflexive Banach spaceand let $A$, $B$ be in $L(X)$.  Then
\begin{itemize}
\item[(a)] If $ATB$ is a compact operator on $X$ for every $T$ in $L(X)$, then $L_A R_B$ is a weakly compact operator on $L(X)$.
\item[(b)] If $X$ has the compact approximation property and $L_A R_B$ is a weakly compact operator on $L(X)$, then $ATB$ is a compact operator on $X$ for every $T$ in $L(X)$.
\end{itemize}
\end{thm}

\pf  To prove (a), it is enough to recall \cite{kalton} that for a reflexive space $X$, on bounded subsets of $K(X)$ the weak topology is the same as the weak operator topology (the operator $T\mapsto f_T\in C((B_X, \mbox{ weak}) \times (B_{X^*}, \mbox{weak}))$, where $f_T(x,x^*):=\langle x^*,Tx\rangle$,  is an isometric isomorphism from $K(X)$ into a space of continuous functions on a compact Hausdorff space).

To prove (b), suppose that we have a $T\in L(X)$ with $ATB$   not compact. Then there is a weakly null normalized sequence $\xn$ in $X$ and $\delta>0$ so that for all $n$, $\|ATB x_n\|>\delta$. Since a reflexive space with the compact approximation property also has the compact metric approximation property \cite{chojohnson}, there are $C_n\in K(X)$ with $\|C_n\|< 1+1/n$, $C_n Bx_i=Bx_i$ for $i\le n$.  Since the $C_n$ are compact, for each $n$, $\|C_n B x_m\| \to 0$ as $m\to \infty$.  Thus $A(TC_n)B x_i=ATB x_i$ for $i\le n$ and $\|A(TC_n)B x_m\|\to 0$ as $m\to \infty$.  This implies that no convex combination of $\{A(TC_n)B\}_{n=1}^\infty$ can converge in the norm of $L(X)$ and hence $\{A(TC_n)B\}_{n=1}^\infty$ has no weakly convergent subsequence.  This contradicts the weak compactness of $L_A R_B$ and completes the proof.
\endpf

\vfill\eject

\vfill\eject
\noindent William B. Johnson\newline
             Department Mathematics\newline
             Texas A\&M University\newline
             College Station, TX, USA\newline
             E-mail: johnson@math.tamu.edu

\bigskip

\noindent Gideon Schechtman\newline Department of
Mathematics\newline Weizmann Institute of Science\newline Rehovot,
Israel\newline E-mail: gideon@weizmann.ac.il


\begin{thebibliography}{99}





\bibitem[Al]{ald} Aldous, David J.  Subspaces of $L\sp{1}$,
via random measures.
 Trans. Amer. Math. Soc.  267  (1981),  no. 2, 445--463.


\bibitem[Bo]{bon} Bonami, Aline. Etude des coefficients de Fourier des fonctions de
 $L\sp{p}(G)$.
(French)  Ann. Inst. Fourier (Grenoble)  20  1970  fasc. 2
335--402 (1971).

\bibitem[CJ]{chojohnson} Cho, Chong-Man; Johnson, William B.
A characterization of subspaces $X$ of $l\sb p$ for which $K(X)$ is an $M$-ideal in $L(X)$.
Proc. Amer. Math. Soc. 93 (1985), no. 3, 466--470.

\bibitem[DJT]{djt} Diestel, Joe;  Jarchow, Hans;  Tonge,
Andrew. Absolutely summing operators. Cambridge Studies in
Advanced Mathematics, 43. Cambridge University Press, Cambridge,
1995.

\bibitem[JJ]{jj} Johnson, W. B.;  Jones, L.  Every
$L\sb{p}$ operator is an $L\sb{2}$ operator.
 Proc. Amer. Math. Soc.  72  (1978),  no. 2, 309--312.

\bibitem[Jo]{j} Johnson, William B.  Operators into $L\sb{p}$
which factor through $\ell\sb{p}$.
 J. London Math. Soc. (2)  14  (1976),  no. 2, 333--339.


\bibitem[KP]{kp} Kadec, M. I.;  Pe\l czy\'nski, A.  Bases,
lacunary sequences and complemented subspaces in the spaces
 $L\sb{p}$.
 Studia Math.  21  1961/1962 161--176.

\bibitem[Kal]{kalton}  Kalton, N. J.
Spaces of compact operators.
Math. Ann. 208 (1974), 267--278.



\bibitem[KK]{kk} Koldobsky, Alexander;  K\"onig, Hermann.
Aspects of the isometric theory of Banach spaces.
 Handbook of the geometry of Banach spaces, Vol. I,
 899--939, North-Holland, Amsterdam,  2001.

\bibitem[Kr]{kr} Krivine, J. L.  Sous-espaces de dimension
finie des espaces de Banach
 reticules.
 Ann. of Math. (2)  104  (1976),  no. 1, 1--29.

\bibitem[KM]{km} Krivine, J.-L.;  Maurey, B.  Espaces de
Banach stables. (French)  [Stable Banach spaces]  Israel J. Math.
39  (1981),  no. 4, 273--295.

\bibitem[La]{la} Lacey, H. Elton. The isometric theory of
classical Banach spaces. Die Grundlehren der mathematischen
Wissenschaften, Band 208. Springer-Verlag, New York-Heidelberg,
1974.

\bibitem[LT]{lt} Lindenstrauss, Joram;  Tzafriri, Lior. Classical Banach spaces
I\&II.  Ergebnisse der Mathematik und ihrer Grenzgebiete, Vol. 92
\& 97. Springer-Verlag, Berlin-New York,  1977 \& 1979

\bibitem[MP]{mp} Maurey, Bernard;  Pisier, Gilles. Series
de variables aleatoires vectorielles independantes et
 proprietes geometriques des espaces de Banach.
(French)  Studia Math.  58  (1976),  no. 1, 45--90.

\bibitem[PR]{pelros}  Pe\l czy\'nski, A.; Rosenthal, H. P.
Localization techniques in $L\sp{p}$ spaces.
Studia Math. 52 (1974/75), 263--289.



\bibitem[Pl]{pl} Plotkin, A. I.  An algebra that is
generated by translation operators, and
 $L\sp{p}$-norms.
(Russian)  Functional analysis, No. 6: Theory of operators in
linear spaces
 (Russian),
 pp. 112--121. Ulyanovsk. Gos. Ped. Inst., Ulyanovsk,  1976.

\bibitem[Ro1]{rosxp} Rosenthal, Haskell P.  On the subspaces of
$L\sp{p}$ $(p>2)$ spanned by sequences of
 independent random variables.
 Israel J. Math.  8  (1970), 273--303.

\bibitem[Ro2]{rosII} Rosenthal, Haskell P.  On the span in
$L\sp{p}$ of sequences of independent random variables.
 II.
 Proceedings of the Sixth Berkeley Symposium on Mathematical Statistics
 and Probability (Univ. California, Berkeley, Calif., 1970/1971), Vol. II:
 Probability theory,
 pp. 149--167. Univ. California Press, Berkeley, Calif.,  1972.

\bibitem[Ro3]{ros} Rosenthal, Haskell P.  On subspaces of
$L\sp{p}$.
 Ann. of Math. (2)  97  (1973), 344--373.

 \bibitem[Ro4]{rosAltgeld} Rosenthal, Haskell P. Convolution by a
 biased coin. The Altgeld Book 1975/76.

\bibitem[Ru]{ru} Rudin, Walter. $L\sp{p}$-isometries and
equimeasurability.
 Indiana Univ. Math. J.  25  (1976),  no. 3, 215--228.

\bibitem[ST1]{st1} Saksman, Eero;  Tylli, Hans-Olav. Weak
compactness of multiplication operators on spaces of bounded
 linear operators.
 Math. Scand.  70  (1992),  no. 1, 91--111.

 \bibitem[ST2]{st2} Saksman, Eero;  Tylli, Hans-Olav.
 Multiplications and elementary operators in the Banach space
 setting. In: Methods in Banach space theory. (Proc. V Conference
 on Banach spaces, Caceres, September 13-18, 2004; J.F.M. Castillo
 and W.B. Johnson, eds.) London Mathematical Society Lecture Note
 Series 337 (Cambridge University Press, 2006), pp. 253-292.

\bibitem[We]{weis} Weis, L.
On perturbations of Fredholm operators in $L\sb{p}(\mu)$-spaces.
Proc. Amer. Math. Soc. 67 (1977), no. 2, 287--292.



\bibitem[Wo]{woj} Wojtaszczyk, P.  Banach spaces for
analysts. Cambridge Studies in Advanced Mathematics, 25. Cambridge
University Press, Cambridge,  1991.


\end{thebibliography}
\end{document}